\documentclass[12pt]{article}
\usepackage{amsmath,amssymb}
\newtheorem{thm}{Theorem}[section]

\newtheorem{cor}{Corollary}[section]
\newtheorem{prop}{Proposition}[section]

\newcommand{\x}{\times}
\newcommand{\bs}{\bigskip}
\newcommand{\dt}{\cdot}

\newcommand{\g}{\gamma}
\newcommand{\G}{\Gamma}

\newcommand{\var}{\varphi}
\newcommand{{\z}}{\Bbb Z}
\renewcommand{\P}{\mathbb P}
\newcommand{\pr}{{\mathbb P}^4}
\renewcommand{\O}{\mathcal O}
\newcommand{\cs}{cubic scroll }
\newcommand{\his}{Hilbert scheme }
\newcommand{\bil}{biliaison }
\newcommand{\li}{liaison }
\newcommand{\cu}{curve}
\newcommand{\I}{\mathcal I}
\newcommand{\irr}{irreducible }
\newcommand{\ex}{example}
\newcommand{\equ}{equivalent }
\newcommand{\co}{codimension }

\begin{document}
\setlength{\baselineskip}{16pt}

\begin{center}
{\Large\bf Experiments with Gorenstein Liaison}

\bs\bs {\sc Robin Hartshorne}

\bs Department of Mathematics\\University of California\\
Berkeley, California 94720--3840

\bs\bs {\it Dedicated to Silvio Greco on his 60th birthday}
\end{center}

\bs\bs
\begin{quote}{\bf Abstract.}
We give some experimental data of Gorenstein liaison, working with
points in ${\P}^3$ and curves in $\pr$, to see how
far the familiar situation of liaison, biliaison, and Rao modules
of \cu s in ${\P}^3$ will extend to subvarieties of
\co 3 in higher ${\P}^4$.\\ \mbox{}\hfill\begin{tabular}{r}
{\sc AMS classification numbers:} 14H50,~14M07\end{tabular}
\end{quote}

\section{The Problem}

For \cu s in projective three-space ${\P}^3_k$, the usual theory of
\li and \bil is well understood \cite{MDP}. We will recall some of the basic
facts, and then explore to what extend these results may generalize to \li
classes of varieties of higher codimension, such as \cu s in ${\P}^4$. Our
method is to run experiments in various special cases, and look for \ex s
which may indicate how the general situation will be.

First we recall the situation in ${\P}^3_k$. A {\it \cu} will be a pure
one-dimensional locally Cohen-Macaulay closed subscheme of ${\P}^3$.
Two \cu s $C_1$ and $C_2$ are {\it linked} if there exists a complete
intersection \cu $D$ such that $D=C_1\cup C_2$ set-theoretically, and
\begin{eqnarray*}
{\I}_{C_1,D} &\cong & {\mathcal H}\!\,om({\O}_{C_2},{\O}_D) \\
{\I}_{C_2,D} &\cong & {\mathcal H}\!\,om({\O}_{C_1},{\O}_D) \ .
\end{eqnarray*}
The equivalence relation generated by linkage is called {\it \li .}
An even number of linkages generates the equivalence relation of
{\it even \li} or {\it \bil.}

We say that $C_2$ is obtained from $C_1$ by an {\it elementary \bil of height}
$h$ if there is a surface $S$ in ${\P}^3$ containing $C_1$, and
$C_2\sim C_1 +hH$ on $S$, where $\sim$ denotes linear equivalence,
and $H$ is the hyperplane section. Here we use the theory of
generalized divisors \cite{GD}, so that any curve on any surface  in
${\P}^3$ can be regarded as a divisor.

Then it is known that an elementary \bil is an even \li\!\!, and the 
equivalence
relation generated by elementary \bil\!\!s is the same as even \li
\cite[4.4]{GD}.

Some of the main results of \li theory for curves in ${\P}^3$ are contained
in the following theorem.

\begin{thm} For curves in ${\P}^3_k$, we have

{\rm{a)}} Two curves $C_1,C_2$ are in the same \li equivalence class
if and only if their Rao modules $M(C_i)=H^1_*({\I}_{C_i}(n))$ are
isomorphic, up to dualizing and shifting degrees. They are in the same
\bil equivalence class if and only if $M(C_1)$ and $M(C_2)$ are
isomorphic up to  shift of degrees.

{\rm{b)}} For each finite-length graded module $M_0$ over the
homogeneous coordinate ring\linebreak $R=k[x_0,x_1,x_2,x_3]$, there exists a
smooth
\irr \cu  \ $C$ in ${\P}^3$ and an integer $h$, such that $M(C)\cong M_0(h)$.

{\rm{c)}} For any finite length $M_0\neq 0$, there is a minimum $h$ for
which there exist \cu s $C_0$ with $M(C_0)=M_0(h)$. These are called
{\bf minimal} \cu s, and the family ${\mathcal L}_0(M_0)$ of universal
\cu s for $M_0$ is an \irr subset of the \his.

{\rm{d)}} (The Lazarsfeld-Rao property): Any other \cu \ $C$ in the \bil class
associated to the module $M_0$ can be obtained by a finite number
of ascending (i.e. $h\geq 0$) elementary \bil\!\!s, plus a deformation,
from a universal \cu \ $C_0$ in the \bil class.

{\rm{e)}} For any module $M$ and any postulation character $\g$, the
subset $H_{\g,M}$ of the \his of \cu s with postulation character $\g$
and Rao module $M$ is \irr (provided it is non-empty).
(For a \cu \ $C$ with homogeneous coordinate ring $R(C)=R/I_C$,
we define the {\it postulation character} $\g_C$ to be the third
difference function of the negative of the Hilbert function
$\var(\ell)=\dim_k R(C)_{\ell}$ of $C$.)
\end{thm}

For proofs of these results, see \cite{R} for a), b), and
\cite{MDP} for c), d), and e).

\bs
A curve $C$ is {\it arithmetically Cohen-Macaulay} (ACM) if its homogeneous
coordinate ring $R(C)$ is a Cohen-Macaulay ring. The ACM \cu s form a special
case of the above theorem that requires slightly modified statements.

\begin{thm}
{\rm{a)}} A \cu \ $C$ is ACM if and only if its Rao module is 0.
The ACM \cu s form one \bil equivalence class.

{\rm{b)}} Any ACM \cu \ can be obtained from a line by a finite
number of ascending elementary \bil\!\!s, plus a deformation.

{\rm{c)}} The postulation character $\g$ of an ACM curve is {\bf positive}
in the following sense: $\g(0)=-1$; if $s_0$ is the least integer $\geq 1$
for which $\g(s_0)\geq 0$, then $\g(n)\geq 0$ for all $n\geq s_0$.
Conversely, for every positive postulation character, there exists an
ACM \cu \ with that character.

{\rm{d)}} If the ACM \cu \ $C$ is integral,
then its postulation character is {\bf connected}, meaning that
$\{n\in{\z}\mid \g(n)>0\}$ is an interval in ${\z}$. Conversely,
for every connected positive character, there is a smooth \irr ACM \cu \
with that character.

{\rm{e)}} For any positive $\g$, the set of ACM \cu s with
postulation character $\g$ is an \irr subset of the \his.
\end{thm}

For proofs, see \cite{E}, \cite{GP}, and \cite{MDP}.

\bs
Now our problem is to what extent do these results extend to \cu s in ${\P}^4$,
or more generally to subschemes of codimension $\geq 3$ in
any projective space?

First of all, it is clear that the definition of \li given above using
complete intersections (which we denote by CI-\li\!\!) is too restrictive.
This has been made abundantly clear in the work of \cite{KMMNP} ---
see the report of R.~Mir\'o-Roig in this volume \cite{Miro}:
there are other invariants besides the Rao module for CI-\li in codimension 3,
and using these, one can  construct many examples of \cu s in ${\P}^4$
having  the same Rao module, but not in the same CI-\li class.

Therefore we will take Gorenstein \li to be the natural generalization of
CI-\li to higher codimension. We state the definitions for \cu s in ${\P}^4$,
though the generalization to subschemes of any dimension in any ${\P}^4$ is
obvious \cite{M}.

A \cu \ $D$ in ${\P}^4$ is {\it arithmetically Gorenstein} (AG) if its
homogeneous coordinate ring $R(D)=R/I_D$ is a Gorenstein ring, where now
$R=k[x_0,x_1,x_2,x_3,x_4]$ is the coordinate ring of ${\P}^4$. Two
curves $C_1,C_2$ in ${\P}^4$ are $G${\it -linked} if there exists an AG
curve $D$ satisfying the same conditions as in the definition of \li  for \cu s
in ${\P}^3$ above. The equivalence relation generated by $G$-linkage is
$G${\it -\li.} The equivalence relation generated by even numbers of
$G$-linkages is {\it even $G$-\li.}

A \cu $C_2$ is obtained from $C_1$ by an {\it elementary G-\bil of height} $h$
if there exists an ACM surface $X$ in $\pr$ satisfying also $G_1$
(Gorenstein in codimension one), containing $C_1$, such that $C_2\sim C_1
+hH$ on $X$, where again $H$ is the hypersurface section of $X$.

It is easy to see that a $G$-\bil is an even $G$-\li\!\! \cite[$\S$5.4]{M}.
The authors of
\cite{KMMNP} are fond of speaking of $G$-\li ``as a theory of divisors
on arithmetically Cohen-Macaulay schemes," and indeed, most of their
\ex s of $G$-\li can also be accomplished by elementary $G$-\bil\!\!s.
However, the relation between these two notions is not yet clear,
so we pose it as a question.

\begin{quote}{\sc Question 1.3} \ \ Is the equivalence relation generated
by elementary
$G$-biliaisons equivalent to even $G$-\li\!\!?
\end{quote}
This is true for CI-\li in any codimension \cite[4.4]{GD},
hence for $G$-\li in codimension 2, but is already unknown for
\cu s in $\pr$.

It is easy to see that evenly $G$-linked \cu s have the same Rao module,
up to twist \cite[5.3.3]{M}, but the converse is unknown:
\begin{quote}{\sc Question 1.4} \ \ If two \cu s $C_1,C_2$ in $\pr$ have
isomorphic Rao modules, up to shift in degrees, are they in the
same $G$-\li class?
In particular, are any two ACM \cu s in the same \bil class?
\end{quote}
(This is \equ to asking if every ACM \cu \ is {\it glicci}, an acronym for
``Gorenstein \li class of a complete intersection.")

For a given finite-length graded module $M\neq 0$, it is easy to see
there is a minimum twist $M(h_0)$ for which there are \cu s with
Rao module $M(h_0)$ \cite[1.2.8]{M}. These are called {\it minimal \cu s.}
Migliore has observed \cite[5.4.8]{M} that the set of minimal \cu s for a given
$M$ may not be \irr\!\!, so we state
\begin{quote} {\sc Problem 1.5} \ \ For a given module $M\neq 0$,
describe the set of minimal \cu s for the module $M$.
Are they all in the same even $G$-\li class?
\end{quote}
As an analogue of the Lazarsfeld-Rao property, we ask
\begin{quote} {\sc Question 1.6} \ \ If $C$ is a \cu \ with Rao
module $M\neq 0$, can $C$ be obtained by a finite number of ascending
elementary $G$-\bil\!\!s from a minimal \cu \ for the module $M$?
For ACM curves we ask, can any ACM \cu \ be obtained by a finite number
of ascending
elementary $G$-\bil\!\!s from a line?
\end{quote}
And lastly,
\begin{quote} {\sc Question 1.7} \ \ Does the set of \cu s with given
Rao module $M$ and postulation character $\g$ form an \irr subset of the
\his?
\end{quote}
In spite of the optimism of some of the researchers mentioned in the 
references,
my expectation is that many of these questions will have negative answers.
The purpose of this talk is to give negative answers to a couple of these
questions, and to propose potential counter\ex s to some others.
We refer to the paper \cite{RH1} for more details of results only stated
here, and further references.

\section{Points in ${\P}^3$}

Closed subschemes of dimension zero of ${\P}^3$ form the first non-trivial
case of \co 3 schemes in a ${\P}^4$. Any such scheme is ACM, so the
questions to consider are a) is every such scheme glicci? and b)
can any such scheme be obtained from a single point by a sequence
of ascending $G$-\bil\!\!s (or ACM curves in ${\P}^3$)?

Since the structure of arbitrary zero-dimensional subschemes can be
quite complicated (unlike the case of zero-schemes in ${\P}^2$, the \his
of zero-schemes of degree $d$ in ${\P}^3$ for fixed $d$ may not even
be \irr\!! \cite{I}), we decided to consider only reduced zero-schemes,
i.e., finite sets of points, in general position. Here {\it general position}
will always mean for a suitable Zariski-open subset of the \his\!\!, possibly
subject to the condition of lying in a given \cu \ or a given surface.
We begin by studying points on low degree surfaces. It is easy to show

\begin{prop} Any set of $n$ general points in ${\P}^2$ can be
obtained by a finite set of\linebreak ascending \bil\!\!s (in this case
CI-\bil is \equ to $G$-\bil\!\!) from a point. {\rm{\cite[2.1]{RH1}}}
\end{prop}

Similarly, using the ACM \cu s on a nonsingular quadric surface,
one can show

\begin{prop} Any set of $n$ general points on a (fixed)
nonsingular quadric surface\linebreak $Q\subseteq {\P}^3$ can be obtained from
  a single point by a finite number of ascending biliaisons
(by ACM \cu s on $Q$). {\rm{\cite[2.2]{RH1}}}
\end{prop}

On a nonsingular cubic surface the situation is more complicated.

\begin{prop} A set of $n$ general points on a (fixed)
nonsingular cubic surface $X\subseteq {\P}^3$
can be connected by $G$-\li\!\!s through sets of general points
of other degrees on $X$ to a single point. In particular a
set of $n$ general points on $X$ is glicci {\rm{\cite[2.4]{RH1}}}
\end{prop}

However, in the proof, we were not able to accomplish this using
ascending \bil\!\!s only. We had to use ascending and descending \li\!\!s
and \bil\!\!s. For \ex\!, to treat 18 general points, one has to link up
to 20, then 28 points, before linking down in many steps to a single point.
\setcounter{cor}{3}
\begin{cor} Any set of $n\leq 19$ general points in ${\P}^3$
is glicci.
\end{cor}

{\bf Proof.} Indeed,  $n\leq 19$ general points  lie on a
nonsingular cubic surface ${\P}^3$.

\bs Our experience in these results is that points lying on surfaces
of low degree 1, 2, or 3, are manageable, but these methods fail
for sets of points on higher degree surfaces. This is consistent with
the \ex s of ACM \cu s in $\pr$, proved to be glicci by \cite[$\S$8]{KMMNP}:
  they lie on rational ACM surfaces which are all contained in
hypersurfaces of degree 1, 2, or 3. So we propose a problem for the
first case not falling under the above results.

\begin{quote} {\sc Problem 2.5} \ \ If $Z$ is a set of 20 points in
general position in ${\P}^3$, is $Z$ glicci? Can $Z$ be obtained by
ascending $G$-\bil\!\!s from a point?
\end{quote}

\section{ACM \cu s in $\pr$}

Following the principle of the previous section, we focus our attention
on general \cu s, usually integral or nonsingular, and sufficiently
general in their component of the \his. Using elementary geometry
of \cu s on the cubic scroll, the Del Pezzo surface of degree 4,
and the Castelnuovo surface of degree 5, we find

\begin{prop} For each possible degree $d$ and genus $g$ of a
nondegenerate integral ACM \cu \ in $\pr$ of degree $d\leq 9$, the
\his $H_{d,g}$ of such \cu s is \irr, and a general such \cu \ can
be obtained by ascending $G$-\bil\!\!s from a line. In particular,
these \cu s are glicci. {\rm{\cite[3.4]{RH1}}}
\end{prop}

A similar argument, using \cu s on the Bordiga surface of degree 6,
gives the same result for ACM \cu s with $(d,g)\!=\!(10,6)$.

\bs{\bf Example 3.2} \ \ For  $(d,g)\!=\!(10,9)$ the \his of smooth ACM \cu s
in $\pr$ has two \irr components. To see this, first consider a
nondegenerate smooth (10,9) \cu \ $C$ in $\pr$. Since $h^0({\O}_C(2))=12$,
we find $h^0({\I}_C(2))\geq 3$. It follows that $C$ is contained in an \irr
surface of degree 3, which must be either a cubic scroll or the cone over a
twisted cubic \cu \ in ${\P}^3$.

We represent the cubic scroll $S$ as ${\P}^2$ with one point blown up.
If $\ell$ is the total transform of a line in ${\P}^2$, and $e$ is the class
of the exceptional line, we denote a divisor $D=a\ell-be$ by $(a;b)$.
Then $S$ is embedded in $\pr$ by $H\!=\!(2;1)$. In this notation there
are two types of smooth (10,9) \cu s, $C_1=\!(6;2)$ and $C_2=\!(7;4)$. Note
that each of these is obtained by $G$-\bil from a line on $S$: $C_1\sim
L_1+3H$ where $L_1=(0;-1)$ and $C_2\sim L_2+3H$ where $L_2=(1;1)$. Hence
both types are ACM.

The two types are distinguished by the following properties:
\begin{quote}
\begin{description}
\item{a)} their self-intersection on  $S$: $C^2_1=32$ while $C_2^2=33$.
\item{b)} their trisecants: since $S$ is an intersection of quadric
hypersurfaces, any trisecant to $C_i$ must lie in $S$. The lines in $S$ are
of types $L_1,L_2$ above. So we see that $C_1$ has no trisecants, while
$C_2$ has infinitely many trisecants of type $L_2$.
\item{c)} their gonality: $C_2$ is trigonal (a $g^1_3$ is cut out by the
trisecants) while $C_1$ is not trigonal.
\item{d)} their multisecant planes. Let $\pi$ be a plane containing
the conic $\G$ of type (1;0) on $S$. Then $C_1\dt\pi=6$ and $C_2\dt\pi=7$.
The pencil of hyperplanes through $\pi$ cuts out a $g^1_4$ on $C_1$
and a $g^1_3$ on $C_2$, computing the gonality of each \cu\!.
\end{description}\end{quote}

Because each of these \cu s is contained in a unique cubic scroll,
if $C_t$ is a family of smooth (10,9) \cu s, then it is contained
in a family $S_t$ of cubic surfaces. Hence the self-intersection of $C_t$
on $S_t$ is constant in a family, and we conclude that neither type can
specialize to the other. Hence the \his of smooth curves $H_{10,9}$
has two \irr components, represented by the two types $C_1$ and $C_2$.

In contrast to the situation in ${\P}^3$ (where for \ex, the \his of
smooth curves of $(d,g)\!=\!(9,10)$ has two disconnected components),
our two  components  of $H_{10,9}$ in $\pr$ have a common intersection,
formed by smooth (10,9) \cu s lying on the singular cubic surface $S_0$,
the cone over a twisted cubic \cu \ in ${\P}^3$. In a family $S_t$ of
smooth cubic scrolls, with limit $S_0$, both classes of lines $L_1,L_2$
have as limit a ruling $L_0$ of the cone $S_0$. So the two divisor classes
$C_1,C_2$ both tend to the singular divisor class $L_0+3H$ on $S_0$.
It is easy to see this divisor class on $S_0$ contains smooth \cu s.
Then, imitating the proof of \cite[2.1]{Z}, cf. \cite[1.6]{Z},
one can show that every smooth (10,9) \cu \ on $S_0$ is a limit of
flat families of \cu s of either type $C_1$ or type $C_2$ on cubic scrolls.

Note finally, since ${\O}_C(2)$ is already nonspecial, it is easy to see that
the postulation of all ACM (10,9) \cu s is the same, so we have an \ex \ where
the \his of ACM \cu s with a fixed postulation is not \irr\!, answering
Question 1.7 above.

To show that this \ex \ is not an isolated phenomenon, we prove the following.

\setcounter{thm}{2}
\begin{thm} Let $X$ be a smooth ACM surface in $\pr$, let $C_0\subseteq X$
be a \cu, and assume either {\rm{a)}} $X$ is rational, or {\rm{b)}}
$C_0\sim aH+bK$ for $a,b\in{\z}$, where $H$ is the hyperplane class, and $K$
the canonical divisor on $X$. Then for $m>>0$, the set of \cu s $C\sim
C_0+mH$ on $X$, together with their deformations $C_t\subseteq X_t$ as
$X$ moves in the family of ACM surfaces $X_t$, forms an open subset of an \irr
component of the \his of \cu s in $\pr$.
\end{thm}

{\bf Proof.} For $m>>0$, each such \cu \ $C$ will lie on a unique such $X_t$,
so the dimension of the family of these \cu s will be equal to the
dimension of the linear system $|C|$ on $X$, which is equal to
$h^0({\mathcal N}_{C/X})$, where ${\mathcal N}$ denotes normal bundle,
plus the dimension of the family of ACM surfaces $X$, which is equal to
$h^0({\mathcal N}_{X/{\P}^4})$ by \cite{E}. (Here the hypothesis
a) or b) of the statement guarantees that when we deform $X$, the
divisor class $C$ extends to the deformed surface.) On the other hand,
we know that the dimension of the family of these \cu s is
$\leq h^0({\mathcal N}_{C/{\P}^4})$ by the differential study of the \his.

Now from the exact sequence
\[
0 \ \to \ {\mathcal N}_{C/X} \ \to \ {\mathcal N}_{C/{\P}^4} \ \to \ {\mathcal
N}_{X/{\P}^4}
\otimes {\O}_C \ \to \ 0
\]
we find
\[
h^0({\mathcal N}_{C/{\P}^4}) \leq h^0({\mathcal N}_{C/X}) +
h^0({\mathcal N}_{X/{\P}^4}
\otimes {\O}_C) \ .
\]
On the other hand, consider the exact sequence
\[
0 \ \to \ {\mathcal N}_{X/{\P}^4}(-C) \ \to \ {\mathcal N}_{X/{\P}^4} \ \to \
  {\mathcal N}_{X/{\P}^4}
\otimes {\O}_C \ \to \ 0 \ .
\]
Since $C\sim C_0+mH$, it follows from duality and Serre vanishing that
$h^i({\mathcal N}_{X/{\P}^4}(-C))=0$ for $i=0,1$ and for $m>>0$.
Hence $h^0({\mathcal N}_{X/{\P}^4})=h^0{\mathcal N}_{X/{\P}^4}\otimes{\O}_C)$
for $m>>0$.

Putting these inequalities together, we find that $h^0({\mathcal 
N}_{X/{\P}^4})$
is equal to the dimension of the family of \cu s in question. We conclude
from this that they form an open subset of a (generically reduced)
\irr component of the \his.

\bs{\bf Example 3.4} \ \ We can use this theorem to make more
examples of non-\irr \his\!\!s of \cu s with given postulation and Rao module.

A first \ex \ is furnished by the families $C_1+mH$ and $C_2+mH$ on the cubic
scroll, where $C_1,C_2$ are the \cu s of Example 3.2.
For given $m$, they will have the same degree, genus, and postulation; each
forms an open set of an \irr component of the \his\!, but the two families
are different because the \cu s have a different self-intersection on $S$.

For another \ex , let $X$ be a Bordiga surface, represented as ${\P}^2$
with 10 points $P_1,\dots ,P_{10}$ blown up, where the notation
$(a;b_1,\dots ,b_{10})$ represents the divisor $a\ell-\sum b_ie_i$,
and the embedding is given by $H=\!(4;1^{10})$. Consider the divisors
\begin{eqnarray*}
L_1 &=& (0;0^9,-1) \\
L_2 &=& (1;1^3,0^7) \\
L_3 &=& (2;1^7,0^3) \ .
\end{eqnarray*}
On a general Bordiga surface, $L_1$ is a line, and $L_2,L_3$ are not
effective. But if $P_1,P_2,P_3$ are collinear, we get a special smooth
Bordiga surface on which $L_2$ is represented by a line. If $P_1,\dots ,P_7$
lie on a conic, we get another special Bordiga surface on which
$L_2$ is represented by a line. It follows that $C_i\sim L_i+mH$ are ACM
\cu s with the same postulation on a general Bordiga surface, for $i=1,2,3$,
and $m>>0$.

By the theorem, each of these $C_i$ forms an (open set of) an \irr component of
the \his\!. Since $L^2_i=-1$, \ $L^2_2=-2$, \
$L^2_3=-3$, the $C_i$ have different self-intersection, so the
components are distinct.

\begin{quote}{\sc Problem 3.5} \ \ Find a way to distinguish the \irr 
components
of the \his of ACM \cu s in $\pr$ with given degree, genus, and postulation.
For \ex\!, would any of the properties suggested in a),b),c),d) of Example 3.2
force the family to be \irr\!?\end{quote}

{\bf Example 3.6} \ \ Our last experiment with ACM \cu s is the first
case of an ACM \cu \ not contained in a cubic hypersurface, namely
smooth ACM \cu s with $(d,g)\!=\!(20,26)$.

There are such \cu s defined by the $4\x 4$ minors of a $4\x 6$ matrix
of general linear forms. These determinantal \cu s are glicci, by a
theorem of \cite{KMMNP} and form an \irr family of dimension $\leq 69$
  \cite[10.3]{KMMNP}.

Allowing these determinantal \cu s  to move in linear systems on smooth ACM
surfaces $X$ of degree 10 and sectional genus 11, we get a family of \cu s of
dimension $\leq 74$, whose general member can be obtained by ascending
$G$-\bil\!\!s from a line, and hence is glicci \cite[3.9]{RH1}.

On the other hand, the differential study of the \his shows that every \irr
component of smooth \cu s of $(d,g)\!=\!(20,26)$ must have dimension
$\geq \!5d+\!1\!-\!g=75$.

By a subtle study of the dimensions of linear systems of \cu s on the
ACM surface of degree 10 mentioned above, we show that a general \cu \
in the \his of (20,26) \cu s cannot be obtained by  ascending
$G$-\bil\!\!s from a line \cite[3.9]{RH1}.
This gives a negative answer to the second half of Question 1.6.
What remains is a problem.

\begin{quote}{\sc Problem 3.7} \ \ Is an ACM \cu \
with $(d,g)\!=\!(20,26)$ in $\pr$ glicci?\end{quote}

\section{Curves on $\pr$ with Rao module $M\neq 0$}

Here the questions to investigate are whether all \cu s with Rao module $M$
belong to the same $G$-\li class; what do the minimal \cu s look like;
and can an arbitrary \cu \  with Rao module $M$
be obtained  from a minimal \cu \ by ascending Gorenstein \bil\!\!s.
As yet, there is very little experimental evidence for these questions,
but what little there is shows that the situation is quite complicated.

We first consider the case $M=k$, of dimension one in one degree only.
We can describe completely the minimal \cu s in this case, which have
$M=k$ in degree 0.

\begin{prop} For every $d\geq 2$ there are minimal \cu s in $\pr$ with Rao
module $M=k$ in degree 0. For each $d$ these \cu s form an \irr family.
The general member of the family is a disjoint union of a line and a plane
\cu \ of degree $d-1$ in general position in $\pr$. Furthermore, all
of these minimal \cu s are in the same $G$-\li class.
{\rm{\cite[4.1]{RH1}}}
\end{prop}

To begin the study of other \cu s with Rao module $M=k$, we look at
smooth \cu s of low degree and genus.
They exhibit many different behaviors.

\bs{\bf Example 4.2} \ \ Every smooth nondegenerate $(d,g)\!=\!(5,0)$ \cu \
in $\pr$ lies on a cubic scroll, has $M=k$ in degree 1, and is obtained by
$G$-\bil from a minimal \cu \ of degree 2, namely two skew lines
\cite[4.3]{RH1}.

\bs{\bf Example 4.3} \ \ Smooth nondegenerate (6,1) \cu s in $\pr$ form an
\irr family. They have $M=k$ in degree 1. They fall into two types. The general
\cu \ $C_1$ lies on a Del Pezzo surface, and is obtained by a $G$-\bil from two
skew lines. This \cu \ has two trisecants.  The special \cu \ $C_2$ lies
on a \cs\!, and is obtained by $G$-\bil from a minimal \cu \ of degree 3.
It has infinitely many trisecants. So in this case the two types are
distinguished by which minimal \cu \ they come from under $G$-\bil
\cite[4.4]{RH1}.

\bs{\bf Example 4.4} \ \ Smooth nondegenerate (7,2) \cu s form an \irr family,
whose general member has $M=k$ in degree 1. In this case the general member of
the family can be obtained by two different routes from minimal \cu s:
one route is $G$-\bil on the Del Pezzo surface from a minimal \cu \ of degree
3; the other is a $G$-\bil on the Castelnuovo surface from a minimal \cu \ of
degree 2 \cite[4.5]{RH1}.

\bs{\bf Example 4.5} \ \ Next we consider smooth nondegenerate (11,7) \cu s
in $\pr$. They form an \irr family,
whose general member has $M=k$ in degree 2. There are such \cu s on a
Bordiga surface, obtained by $G$-\bil in two steps: from two skew lines
to a smooth (5,0) \cu \ on a \cs, then to the (11,7) \cu \ on the Bordiga
surface. However, we can show by counting dimensions that the
general (11,7)  \cu \ does not lie on a Bordiga surface, and cannot be
obtained by ascending $G$-\bil\!\!s from a minimal \cu. This provides
a negative answer to the first part of Question 1.6 above \cite[4.7]{RH1}.
There remains a problem.

\begin{quote}{\sc Problem 4.6} \ \
Is a general (11,7) \cu \ in $\pr$ in the $G$-\li class of
two skew lines?\end{quote}

Minimal \cu s in $\pr$ with Rao module $M_a=R/(x_0,x_1,x_2,x_3,x_4^a)$
for $a\geq 2$ have been studied by Lesperance \cite{L}. He shows

\setcounter{prop}{6}
\begin{prop}
For $a\geq 2$, there are minimal \cu s with Rao module $M_a$ of every
degree $d\geq a+1$. A reduced minimal \cu \ is one of the following
(where we denote by $P$ the point $(0,0,0,0,1)$.)

{\rm{a)}} A disjoint union of a line and a plane \cu \ of degree $a$
in ${\P}^3$, where $P$ is the point of intersection of the line and the plane.

{\rm{b)}} A disjoint union of  plane \cu s  of degrees $a,b$,
with $a\leq b$, where $P$ is the point of intersection of the two
planes, and $P$ does not lie on either \cu.

{\rm{c)}} A disjoint union of  plane \cu s of degrees $a,b$, with $b\geq 1$,
  where $P$ is the point of intersection of the two planes, but this time
$P$ lies on the \cu \ of degree $b$. (For $b=1$, we recover type {\rm{a)}}
above.)

{\rm{d)}} A disjoint union of a line and an ACM \cu \
in ${\P}^3$, where $P$ is the point of intersection of the line and the
${\P}^3$, and $a$ is the least degree of a surface in ${\P}^3$ containing
the ACM \cu, but not containing $P$.
\end{prop}

{\bf Example 4.8} \ \ In particular, the set of minimal \cu s of a given
degree may not be \irr\!. The first \ex \ is $a=2$, degree 4, where there
are minimal \cu s of type b), a union of two conics, and type d), a
union of a line and a twisted cubic \cu, which form two \irr families
\cite[4.5]{L}.

\bs A more serious problem arises with the question of $G$-\li\!.
Lesperance is able to show that most of the minimal \cu s described
in Proposition 4.7 are in the same $G$-\li class as the first (type a).
However there remains an open question, of which we state the first case.

\begin{quote}{\sc Problem 4.9} \ \ Let $C_1$ be a disjoint union of
two conics of type b) above, and let $C_2$ be a disjoint union of a line
and a twisted cubic \cu, of type d) above. Then both have Rao module $M_2$.
Are they in the same $G$-\li class?
\end{quote}

{\bf Example 4.10} \ \ Applying $G$-\bil on a Del Pezzo surface, we can
rephrase Problem 4.9 in terms of smooth \cu s with $(d,g)\!=\!(8,3)$.

On the Del Pezzo surface $X$, note that the divisor class $(1;1,0^4)$ is a
conic, and two such are disjoint. So we can take $C_1=(2;2,0^4)$ on $X$,
and let $D_1=C_1+H=(5;3,1^4)$. This is a smooth (8,3) \cu.

On the other hand, $(1;0^5)$ is a twisted cubic, and $(0;0^4,-1)$ is
a line not meeting it, so we can take $C_2=(1;0^4,-1)$, and $D_2=
C_2+H=(4;1^4,0)$. This is another smooth (8,3) \cu.

The \cu s of type $D_1,D_2$ both have Rao module $M_2$, but neither type
can specialize to the other, because each lies in a unique Del Pezzo surface,
and on that surface, their set of intersection  numbers with the sixteen
lines are $(1^8,3^8)$ for $D_1$ and $(0,1^4,2^6,3^4,4)$ for $D_2$.

Note also that the \his of smooth (8,3) \cu s in $\pr$ is \irr\!,
but the general \cu \ has Rao module $M=k$ in degree 1, and does not
lie on a Del Pezzo surface. Thus our two families of \cu s are locally
closed \irr subsets of $H_{8,3}$.

Both types of \cu s $D_1,D_2$ have self-intersection 12. However,
the two types can be distinguished by
\begin{quote}\begin{description}
\item{a)} their intersections with the 16 lines on $X$ (mentioned above)
\item{b)} their multisecants: $C_1$ has trisecant lines, but no quadrisecant,
while $C_2$ has a quadrisecant line
\item{c)} their multisecant planes: let $\pi$ be the plane containing
the conic $(2;0,1^4)$ in $X$. Then $C_1\dt\pi=6$ while $C_2\dt\pi=5$.
\item{d)} their gonality: $C_1$ is hyperelliptic, with a $g^1_2$ cut out
by the hyperplanes through $\pi$, while $C_2$ has gonality 3, and a
$g^1_3$ is  cut out
by the hyperplanes through $\pi$
\item{e)} the point $P$ (determined by the Rao module) lies on the
surface $X$ for type $D_2$, but does not lie on $X$ for type $D_1$.
\end{description}\end{quote}
Now we can rephrase Problem 4.9 as

\begin{quote}{\sc Problem 4.11} \ \
Do the two types of smooth (8,3) \cu s with Rao module $M_2$
(described above) belong to the same $G$-\li class?
\end{quote}

\section{Conclusion}

For ACM \cu s in $\pr$, we have shown that the family of ACM \cu s with
given degree, genus, and postulation may not be \irr (3.2);
we have given examples of ACM \cu s that cannot be obtained by ascending
Gorenstein \bil from a line (3.6); and we have proposed examples of ACM \cu s
that may not be glicci (3.7).

For \cu s with Rao module $M\neq 0$, we have described the minimal \cu s
in two cases, illustrating their complexity (4.1),(4.7); we have given \ex s
of \cu s that cannot be obtained from a minimal \cu \ by ascending
$G$-\bil\!\!s (4.5); and we have proposed \ex s of \cu s with the same
Rao modules that may not be in the same $G$-\li class (4.6),(4.9).

We have seen by \ex \ that certain families of \cu s with the same Rao module
can be distinguished by the least degree of an ACM surface containing
the \cu, or their self-intersection on an ACM surface of least degree
containing the \cu, or their multisecant lines, or their multisecant planes,
or their gonality. What is lacking at this point is a better understanding of
how these geometrical properties of the \cu \ in its embedding behave under
the operation of Gorenstein \li\!.

\end{document}